\documentstyle[12pt,aps]{revtex}
\begin{document}
\title{The Periodic Response of Periodically Perturbated \\
Stochastic Systems}
\author{A. Yu. Shahverdian}
\address{Yerevan Physics Institute, Yerevan, Armenia}
\author{A. V. Apkarian}
\address{SUNY Upstate Medical University,
Syracuse, NY, USA}
\maketitle
\section{Introduction}
In paper [1] a possible way was noticed regarding how
a neuron retains or responds to
periodic stimulation.
This new feature was established as a result of computational
study of the deep differential structure of actual
neuron spike trains. We observed, that
for many cases the sequences of higher order finite
differences, taken from periodically stimulated
neuron spike trains, can be divided
into several subsequences of approximately equal
lengths, on some of which the changes in monotony of
these differences are strongly periodic. On the other hand, the
presence of chaotic dynamics in neural activity is well known.
This work investigates
the following problem: In which extent the deterministic
stochastic systems hold such type of periodic response property,
or does our remark from [1] remain
valid for living matter only?
With this end in view, we introduce a new numerical
characteristic of stochastic one dimensional systems, expressing
this kind of hidden and incomplete periodicity in a
quantitative form. Its computational study on instances of some
nonlinear systems is presented.
We have chosen three one dimensional maps,
most frequently mentioned in nonlinear dynamics:
the tent map, logistic map, and Poincare displacement of
Chirikov's standard map. The tent map is the simplest
system exhibiting
the strong chaotic properties, while
the behavior of logistic
map is typical for many one dimensional systems having
quadratic maximum. The Chirikov map [2-4]
has a fundamental role in Hamilton dynamics.
We study the behavior of one of its ordinates,
so-called [5] Poincare displacement.
The Billingsley-Eggleston formula type inequality
as well as a
theoretical result on sequences of fractional parts,
generated by a simple aperiodic map, are also presented.
We recall [6], that namely these sequences are well
accepted for the aims of mathematical and computational
modeling of randomness.

Note, that in many problems of nonlinear dynamics
(e.g. various neuron
mathematical models, Duffing equation -- will be studied in
our subsequent work) namely a one dimensional Poincare
displacement is of the main interest.
Let us also note, that the research of various
two dimensional maps and three dimensional flows
(for instance, the orbits of R$\ddot{o}$ssler and Lorentz
attractors - see [2], Ch.~1.5.b and Ch.~7.1.b) can be
mostly reduced to studying some one dimensional systems.
In this regard, we recall Bridges-Rowlands general
theoretical scheme, can be found in [2], Ch.~7.3.b.

We give some statements of finite difference
analysis [1, 7], on which this paper
is based, and explain our approach.
The mentioned new notion,
the numerical characteristic $\gamma$,
is some measure of (in certain sense)
minimal periodicity, and hence,
can also be treated as a measure of chaosity.
The main claim of this work is that
this quantity, being a measure of irregularity, is also
that characteristic of stochastic system, that
is able to change essentially its numerical value when
the system undergoes a
weak periodic perturbation.
This feature of $\gamma$ makes it a new tool when
researching the problems of weak signal detection in presence of
strong (deterministic) noise. In particular, it can be used
in research of stochastic resonance phenomena (see, e.g., [8]).

Accepted approaches in stochastic resonance problems
are the use of power spectrum and correlation function.
However, in this work we give only the comparisions
of new measure with Lyapunov exponent, leaving the consideration
of spectral characteristics for the further work.
The criticism on Lyapunov exponent one can find, e.g., in [9].
Because of its computation simplicity, the $\gamma$
characteristic is well adapted for research of various
applied problems, where the analytic law of system
evolution usually remains unknown.
In this respect, we note
important works by J.~Kurths and his colleagues
([10]; see also [11] where the approach
from [1, 7] to study the seismic time series is applied)
on nontraditional measures and its
applications to studying different chaotic time series in medicine and
astronomy.

\section{The absolute finite differences and a new measure of
irregularity}

The differential method, suggested in [1, 7], reduces the research of
chaotic properties of the orbits
${\bar X}=(x_{i})_{i=1}^{\infty}$ generated by a given one
dimensional system
to analysis of alternations of the monotone increase and 
decrease of higher order absolute finite differences
$$ \Delta^{(s)}x_{i} =
|\Delta^{(s-1)}x_{i+1} - \Delta^{(s-1)}x_{i}|
\qquad (\Delta^{(0)}x_{i}= x_{i} \ ;
\ \ i,s = 1,2,3, \ldots) \ .
$$
In this section we introduce some statements of this approach and
briefly describe some computations, explaining the basic
notions involved in this work.

Let us have an one dimensional system, generating numerical
sequences
${\bar X} = (x_{i})_{i=1}^{\infty}, \ 0\leq x_{i} \leq 1 $
from interval $(0,1)$; we impose no restrictions on
the system, and the generating mechanism can be quite arbitrary.
Following [1, 7] we consider a special representation of finite
orbit ${\bar{X}}_{k}=(x_{i})_{i=1}^{k}$,
which emphasizes its deep differential structure.
Indeed, it can be easly obtained, that
for $1\leq s\leq k-1$ we have
\begin{equation}
  \Delta^{(s-1)}x_{i} = \mu_{k,s-1} +
  \sum_{p=1}^{i-1}(-1)^{\delta_{p}^{(s)}}\Delta^{(s)}x_{p}  -
  \min_{0\leq i\leq k-s}(\sum_{p=1}^{i}
   (-1)^{\delta_{p}^{(s)}}\Delta^{(s)}x_{p})
\end{equation}   
where 
 $$     \delta^{(s)}_{p}= \left \{ \begin{array}{ll}
        0 & \Delta^{(s)}x_{p+1} \geq \Delta^{(s)}x_{p} \\
        1 & \Delta^{(s)}x_{p+1} <    \Delta^{(s)}x_{p}
	                          \end{array} \right.  
    \qquad  \  \mu_{k,s} =   \min \{ \Delta^{(s)}x_{i}:
 1\leq i \leq k-s\} \ ,
 $$
and it is supposed that $\sum_{1}^{0} = 0$.
Hence, one can consider that finite orbits with length
$k$ are given in some special form:
\begin{equation}
{\bar \zeta}_{k} = (r_{1}^{(k)},r_{2}^{(k)},\ldots,r_{m}^{(k)};
   \mu_{k,1}, \  \mu_{k,2}, \  \ldots , \mu_{k,m};
   \rho_{k,1}, \  \rho_{k,2}, \  \ldots, \rho_{k,k-m})
\end{equation}
where 
 $    r_{s}^{(k)} = 0.\delta_{1}^{(s)}\delta_{2}^{(s)}\cdots
                  \delta_{k-s}^{(s)} \quad (1\leq s\leq m)
 $
and
  $$
  \mu_{k,1}, \  \mu_{k,2}, \  \ldots ,
  \mu_{k,m}  \quad \mbox{and} \quad
  \rho_{k,1}, \  \rho_{k,2}, \  \ldots , \rho_{k,k-m}  \qquad
   (\rho_{k,i} = \Delta^{(m)}x_{i})
  $$
are some numbers from interval $[0,1]$. Here $m=m_{k}$ are
some given numbers which tend to $\infty$ as $k\to \infty$
and everywhere below we have chosen $m_{k}=[k/2]$.
It is important to note, that after applying the recurrent
procedure (1) the sequence ${\bar{X}}_{k}$
can be completely recovered by ${\bar{\zeta}}_{k}$. 
For given orbit $\bar X$, we consider
the binary sequences
\begin{equation} 
 {\bar X}^{(n)} = (\delta_{1}^{(n)}, \delta_{2}^{(n)}, \ldots,
 \delta_{k}^{(n)}, \ldots)  \qquad (n, k = 1, 2, \ldots).
\end{equation} 
The method from [7] reduces the study of orbits
$\bar X$ to analysis of some conjugate orbits
${\bar\nu}=(\nu_{n})_{n=1}^{\infty}$ which terms are defined as
follows:
  $$
  \nu_{n} =
  0.\delta_{1}^{(n)}\delta_{2}^{(n)}\delta_{3}^{(n)}\ldots \ \ .
  $$
It is shown [1, 7], that provided some rather general conditions,
the conjugate orbits are attracted to some Cantor set $\cal A$ of
zero Lebesgue measure.

The main claim is that the periodic response of (weak)
periodic perturbation of given stochastic system is localized in
${\bar X}^{(n)}$ and  its presence can be detected,
when considering these sequences.
More exactly, for that purpose it should be studied the
asymptotical (as $N\to \infty$) relative
volume (i.e. the density in natural series) of the set of all
those indeces $i$, for which the changes of binary symbol occur,
\begin{equation}
  \delta_{i+1}^{(N)} = 1-\delta_{i}^{(N)}
\qquad (1\leq i\leq N-1).
\end{equation}
In order to explain how
this statement relates to irregularity
notion and how the transition to chaos occures,
let us consider, from such a
differential point of view, one of our
particular chaotic systems - the logistic map $x\to rx(1-x)$
(it is assumed $0<x<1$ and $0<r<4$).
It is well known, that numerical interval $[0,4]$ is divided into
two subintervals
$[0, r_{\infty})$ and $[r_{\infty},4]$ \
($r_{\infty}=3.569\ldots$),
where the orbits ${\bar X}$ of consecuitive
iterates of this map demonstrate regular periodic
and stochastic and aperiodic motion respectively.
We observed how frequently, in dependence of
the control parameter value $r$, the
changes (4) of binary terms from
Eq.~(3) occur. By this computational way
one can find that for each $r\in [0, r_{\infty})$
there exists some index
$n_{r}$ such that for all $N \geq n_{r}$ the finite
sequence ${\bar X}^{(N)}_{N}$ contains the
series with the same binary
symbol having the lengths tending to $\infty$
as $N\to \infty$. However, when $r$ increasingly
approaches to $r_{\infty}$, these lengths are
decreased and for $r=r_{\infty}$ they become upper
bounded for all $N$. It can be proved [7, 17],
that if they are bounded by a number $K\geq 2$ then
for Hausdorff dimension of attractor $\cal A$ we have
  $$dim ({\cal A}) \leq 1-\frac{1}{(4\ln 2)(K-1)} .$$

Let us now introduce, based on such kind of computational
experience, the following characteristic
of one dimensional systems: for an orbit ${\bar X}=
(x_{k})_{k=0}^{\infty}$ of given system we define
\begin{equation}
 \gamma = \gamma({\bar X}) = {\lim}_{N\to \infty}
 \frac{\gamma({\bar X}, N)}{N}.
\end{equation}
Here, $\gamma ({\bar X}, N)$ denotes the total number of those
indeces $1\leq i\leq N-1$ for each of which equation (4) holds.
Further, we consider the comparisions of $\gamma$ with Lyapunov
exponent $\lambda$ (see e.g., [12], Ch.~5 and [2], Ch.~7.2.b): if
$x_{k+1}=F(x_{k})$, then
\begin{equation}
  \lambda = \lambda({\bar X})=\lim_{N\to
  \infty}\frac{1}{N}\sum_{k=1}^{N}\ln|\frac{dF(x_{k})}{dx_{k}}|.
\end{equation}
This quantity, along with
Kolmogorov-Sinai entropy (that for one dimensional
system coincides with $\lambda$),
power spectrum and correlation function
is one of the most propagated
measures of chaosity.
As the Lyapunov exponent, $\gamma$
exhibits a weak dependence on initial value $x_{0}$.

For the processes, which can be reduced [7] to
that ones generating binary sequences, we have established
(see [17]) the following Billingsley-Eggleston [13] formula type
relation: for Hausdorff dimension of attractor $\cal A$ we have
  $$dim ({\cal A}) \leq -H(\gamma)$$
where $\gamma$ is the system response charactersitic
defined by Eq.~(5), and
  $$H(x)=x\log_{2}x+(1-x)\log_{2}(1-x) \qquad (0<x<1) $$
is the Shannon entropy function.

If $\bar X$ is either constant or periodic, then we obviously have
\begin{equation}
   \gamma({\bar X}, N)=\gamma N+O(1) \qquad (N\to \infty)
\end{equation}
and $0\leq \gamma\leq 1$ is rational. According to next theorem
(see [7, 17]) this also represents a sufficient
condition of regularity in sense of definition from [1] (given
for theoretical neuron spike trains).
\newtheorem{guess}{Theorem}
\begin{guess}
For the sequence ${\bar X} = (\{\alpha n\})_{n=1}^{\infty}$
of fractional parts, where $0<\alpha <1$ is irrational,
the next statements are true:
(a) the conjugate to $\bar X$ orbit
$\bar\nu= (\nu_{n})_{n=1}^{\infty}$ is a periodic
sequence;
(b) if entire part of $1/\alpha$ is of
the form $[1/\alpha]= 2^{p}-1$ ($p\geq 1$)
then $\nu_{n}\equiv 0$ for all enough large indeces $n$.
\end{guess}
This implies, that relation (7) holds for sequences of fractional
parts as well. It is well known [2, 16], that through use
of some canonical transformation, accepted in classical
mechanics,
an arbitrary integrable Hamilton system in fact is
reduced to some billiard system in cube.
In its turn, the billiards boundary behavior, due to
results from [14] (see also [7, 15]), is reduced
to countable set of sequences of fractional parts. It should be also
noticed that the Lyapunov exponent of every integrable system is
equal to zero [2, Ch. 5.3].
In this respect, Theorem 1 and results from [14],
imply that every integrable Hamilton system appears to be
a regular system in some generalized sense of
above mentioned definition from [1] (see [17] for
rigorous formulations).

We note especially, that we are more interested not in
numerical value of quantity $\gamma$ -- our main interest is
focused on the value of its discrepancy when applying to given system
a weak periodic perturbation. The computational
analysis, partially presented in next section, shows that
$\gamma$ possesses the following basic properties:

(A)~~For the most irregular systems $\gamma$ is positive;

(B)~~When applying to irregular system any small
periodic perturbation, $\gamma$ is increased;

(C)~~Different systems have different rate of increase of
$\gamma$.

The second point of Theorem 1 implies that there exist
aperiodic (but integrable) systems for which $\gamma$
is zero. The tent and
logistic maps, after applying a stimulation with  a small
intensity $10^{-4}$ demonstrate an
increase of $\gamma$ up to $25\%$ and $30\%$ respectively, while
for $\theta$-ordinate of standard map the increase
rate can be more than $45\%$.

\section{Computational study of $\gamma$-characteristic}

The introduced above response coefficient $\gamma$ reflects the
asymptotical measure of irregularity of a given
system's orbits in respect of changes in
monotony of finite-differences,
taken from the orbit. In the context of remarks from [1] on
neural activity, it can be said that $\gamma$ also
expresses the
degree of ability of given system "to feel" the stimulation
and to respond on it as the living matter does.
We study the response properties of three parametric
stochastic systems --- the tent map $F_{T}$:
$$F_{T}^{(t)}(x) = t(1-2|\frac{1}{2}-x|) \  \
   (= \left \{ \begin{array}{ll}
   2tx,    & 2x\leq 1\\
   2t-2tx, & 2x>1
   \end{array} \right.)  \qquad (0<x<1; \ 0<t \leq 1),
   $$
the logistic map $F_{L}$:
 $$ F_{L}^{(r)}(x)=rx(1-x) \qquad (0<x<1;\  0<r \leq 4), $$
and the Poincare $\theta$-displacement
$\bar X = (\theta_{n})_{n=1}^{\infty}$ of Chirikov's
standard map $F_{S}^{(K)}=F_{S}^{(K)}(I,\theta)$:
  $$  \begin{array}{ll}
       I_{n+1} & =I_{n}+K\sin \theta_{n}\\
  \theta_{n+1} & = \theta_{n}+I_{n+1} \quad mod\ (2\pi).
  \end{array} \qquad (0<I,\theta<2\pi; \  K>0)
  $$
Here $t$, $r$ and $K$ are control parameters. It is well
known that for $0<t<1/2$ and $0<r<r_{\infty}$ respectively,
the tent and
logistic maps iterates has regular behavior while in
rest part of parameters they mostly exhibit irregular
chaotic motion. The $K>0$ in map $F_{S}$,
so-called stochasticity coefficient, parametrizes
transition from local chaos ($K\approx 0$) to
global "stochastic sea" ($K\approx 1$) [2, 4].

We have been considering also
actual neuron spike trains, obtained from electrophysiological
recordings [18]. However, the data appear to be too
short (note that the work [18] followed other aims)
in order to determine the values of $\gamma$ with a sufficient
accuracy.

Let us now describe the results obtained and to compare the
$\gamma$-characteristic and Lyapunov exponent $\lambda$.
We note, that (see. e.g., [2, 12, 16])
 $$\lambda_{T}(t)=\log_{2} 2t \qquad (0<t \leq 1), \qquad
 \lambda_{S}(K)=\ln\frac{K}{2} \qquad (K\approx 6)$$
while $\lambda_{L}(r)$ has a complex behavior on numerical
interval $3<r<4$; here, the second relation where
$\lambda_{S}$ denotes the maximal Lyapunov exponent of
standard map, has been obtained
by Chirikov for large values of $K$ (see also [12], Ch.5).
After computations we found that different systems may
have different values of
$\gamma$.
The computations show (see Fig.~2) that
for logistic map $F_{L}$ and tent map $F_{T}$ the relation
  $$ \gamma(r)=\alpha(r)\lambda^{+}(r) \qquad
  (\lambda^{+}=\max\{0,\lambda\})$$
holds for all values of control parameter $r$.
Here, $0<\alpha(r)<\infty$ is some continuous function,
which zeroes
can be situated only in zero points of $\lambda^{+}(r)$.
Indeed, one can see (Fig.~2), that
our new measure $0\leq \gamma \leq 1$ always reaches
its local minimal values in small intervals, containing
zero points of $\lambda^{+}$.

More or less exact (with preciseness about $10^{-3}$)
computation of $\gamma$ needs quite large
number of iterates of given map (tens of thousands).
On the other hand, the new measure has that important
advantage,
that is the simplicity of its computation.
It can be easily implemented just over the data
$\bar X$, without refering to process generation law.
In contrary, this
cannot be said about Lyapunov exponent for which, in order to get
an approximate to its numerical value, it is usually
required either a given system's explicit analytic form
(cp. Eq.~(6)) or rather complex theoretical
constructions ([2], Ch.~5.3).
By this reason, as for computation of $\gamma$ we need only
the corresponding time series to be available, the $\gamma$
is better adapted for research of various applied
nonlinear problems.

We were studying the influence of simplest periodic perturbation
on the numerical value of $\gamma$. Namely,
let $(s_{n})_{n=1}^{\infty}$ be a periodic sequence of the form
  $$ s_{n}=  \left \{   \begin{array}{ll}
\epsilon &   n/\tau \ \mbox{is integer} \\
       0 &   n/\tau \ \mbox{is fractional},
       \end{array} \right. $$
i.e. we let $s_{n}=\epsilon$ for the numbers $n$ of the
form $n= \tau, 2\tau, 3\tau, \ldots$ and $s_{n}=0$ elsewhere in
natural series. Here, $0<\epsilon<1$ and natural $\tau \geq 2$
are pregiven. We call such a
sequence $(s_{n})_{n=1}^{\infty}$ the perturbation
(or stimulation) with intensity $\epsilon$ and period $\tau$.
We have been studying the additively perturbated systems
  $$
  {\bar X}_{\epsilon, \tau} =
  (x_{n}^{\prime})_{n=1}^{\infty}
  \quad \mbox{where} \quad
  x_{n+1}^{\prime}= F(x_{n}^{\prime})+s_{n},
  $$
where $F$ denotes one of three mentioned systems.
Such a perturbated system where $F$ is logistic map
but assuming that $s_{n}$ is white
noise, was studied earlier in series of works ([19]; see
details in [16], Ch.~3). In this work we deal with
an inverse statement, considering $F$ as a (deterministic)
noise source and $s_{n}$
as a regular signal. For those actual systems,
for which the analytic shape of generating law
$F$ remains unknown (for
instance, earthquake time series or neuron spike trains;
an explanation of quite complex neuron stimulation
mechanism can be found in [20]), one can
consider ${\bar X}_{\epsilon,\tau}$ as the stimulated system.

Further, for mentioned systems $F$ we have been
studying by computational way the
properties of the response coefficient
  $$ \gamma(\epsilon,\tau)
  = \gamma({\bar X}_{\epsilon,\tau}). $$
First, for different values of intensity $\epsilon$ we have
considered the problem: how large can be the value of
function ${\gamma}(\epsilon, \tau)$
for a given intensity $\epsilon$. Particularly, whether
for given level of intensity $\epsilon$ there exist such
values of period $\tau$, when the response $\gamma$
is close to its possible maximal value $1$?
The result obtained, which demonstrate the
Figs.~3 and 4, can be formulated in the
following form:
\begin{quote}
for given $\epsilon$ the function ${\gamma}(\epsilon, \tau)$ has
the self-affine structure;
\end{quote}
\begin{quote}
for given $\epsilon$ those values of $\tau$,
where the function $\gamma(\epsilon, \tau)$ reaches
its maximal (in respect of $\tau$) value,
are spreaded on the whole natural series and
possess a positive density in this series;
\end{quote}
\begin{quote}
the maximal (in respect of period $\tau$) response
depends on stimulation intensity $\epsilon$ in nonpredictable
way.
\end{quote}
The second statement implies an important conclusion,
that for a given $\epsilon$
the "maximal" period can be found in natural series
with a positive probability.
The Fig.~3 presents the graph of function
$\gamma (\epsilon,\tau)$  for perturbated systems $F$ with applied
stimulations of intensity $\epsilon =0.0001$ and periods
$\tau=2,3,\ldots, 100$, computed with a supercomputer
(we used SiliconGraphics);
the control parameter values are: $r=0.7$ for tent map,
$r=3.7$ for logistic map, $K=0.6$ for standard map.
Compairing graphs on Figs. 1 and 3, one can see that the value of
$\gamma$, which for these three non-perturbated systems are
approximately equal
$0.49$, $0.62$ and $0.35$ (see Fig.~1), under stimulation
is increased (for different periods $\tau$) up to $0.62$, $0.80$,
and $0.50$ respectively.
In this connection, let us note
the following: if to apply to perturbated system $F$
another periodic stimulation with a small
intensity, it should results a similar increase of $\gamma$
(of already once perturbated system). If to continue this process,
it seems possible to construct for a given stochastic
deterministic system a
(multiperiodic) stimulation of any small positive
intensity, which will increase the $\gamma$ up to its
possible maximal value $1$.

We have also studied the maximal coefficient
\begin{equation}
{\gamma}(\epsilon) =
  \max_{2\leq \tau < \infty}{\gamma}(\epsilon,\tau)
\end{equation}
that gives the dependence of maximal values of response on
stimulation intensity $\epsilon$.
This function also possesses
the self-affine shape and it is probably possible to study the
corresponding fractal characteristics.
On the other hand,
one can
see that this function has certain general tendency to monotone
decrease with growth of $\epsilon$. This pecularity
can be emphasized in the following way, accepted in function
theory and classical mechanics: we consider the
averaged function $\mu_{s}$,
 $$ \mu_{s}(\epsilon) =
 \frac{1}{s}\int_{\epsilon}^{\epsilon +s}
 \gamma(t)\delta_{s}(\epsilon,t)dt $$
where
  $$\delta_{s}(x,y)= \left \{ \begin{array}{ll}
          1 & |x-y| < s \\
          0 & |x-y| \geq s
         \end{array} \right. $$
and $s>0$ is some number. Computations show that for enough
small $s>0$ the $\mu_{s}(\epsilon)$
is a monotone decreasing function. This means, that
the maximal (averaged) response of the (above considered)
nonlinear systems is found in the area of small intensities.

\newpage

\section{Figure Captions}

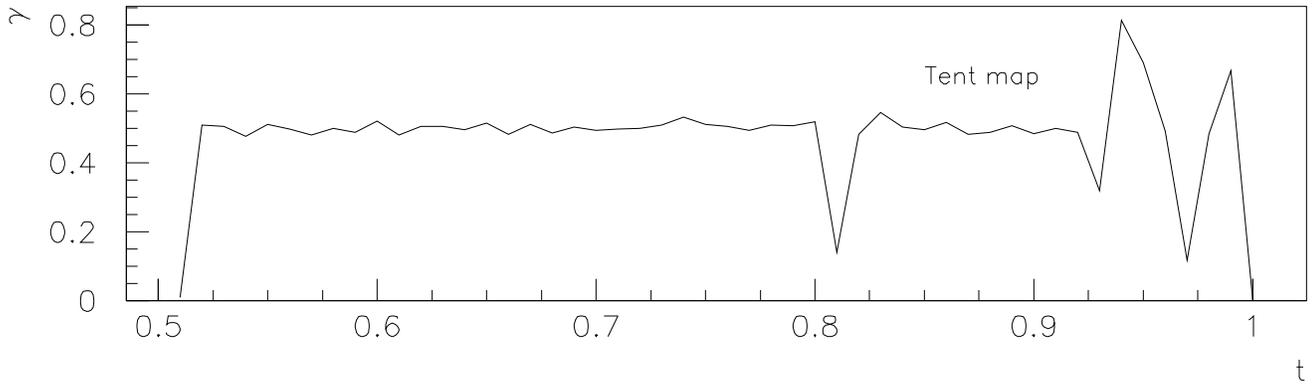
Fig.~1. The graphs of $\gamma$-characteritic
for three different systems $F$,
constructed through relation (5), where is taken
$N=30000$  for tent and logistic maps, and $N=40000$ for
standard map. The value of
$\gamma$ has been computed for 50 values of control parameters
$t$, $r$ and $K$ with the step $0.01$ and $0.02$ respectively.

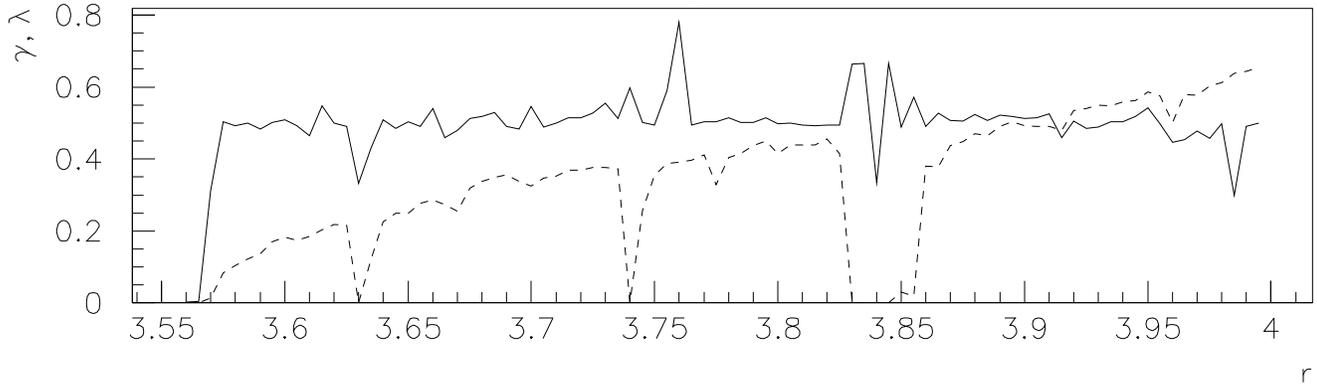
Fig.~2. The graphs of $\gamma$-characteristic
(solid line) and Lyapunov exponent
(dotted line) for logistic map with initial
value $x_{0}=0.55$. (a) Computations made for
$N=50000$ and values of control parameter, starting from
$r=3.56$ till $r=4$ with the constant step is equal $0.005$.
(b) Computations made for $N=50000$, starting from $r=3.7$
till $r=3.9$ with the step $0.002$.
(c) The same computations as in (b) with step 0.001.

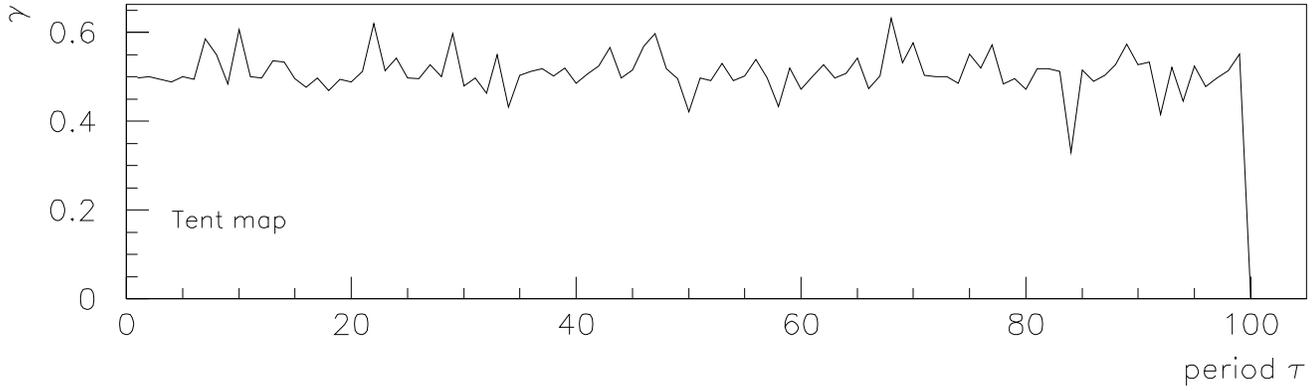
Fig.~3. The graphs of function $\gamma(\epsilon,\tau)$ for three
different systems. It is taken $\epsilon=0.0001$ and
computations made by Eq.~(5)  and for the same
$N$ as in Fig.~1. The systems are: the tent map's iterates
for $t=0.7$ and $x_{0}=0.17$, the logistic map for $r=3.7$ and
$x_{0}=0.317$, and the standard map for $K=0.6$ and
$I_{0}=0.5$, $\theta_{0}=0.2$.
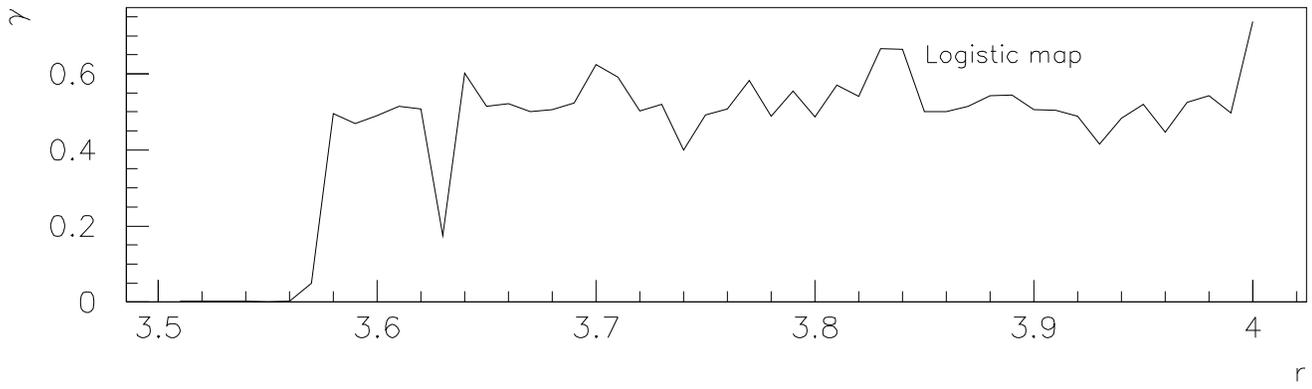
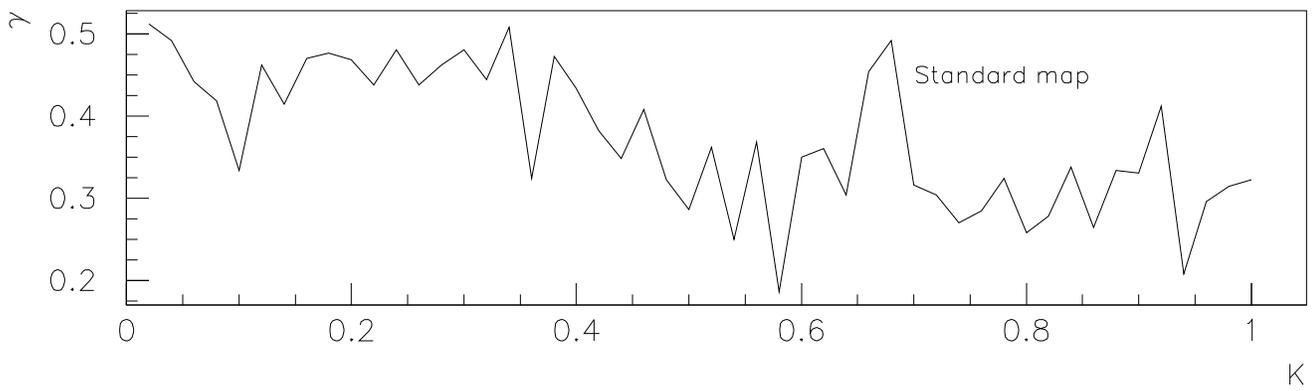
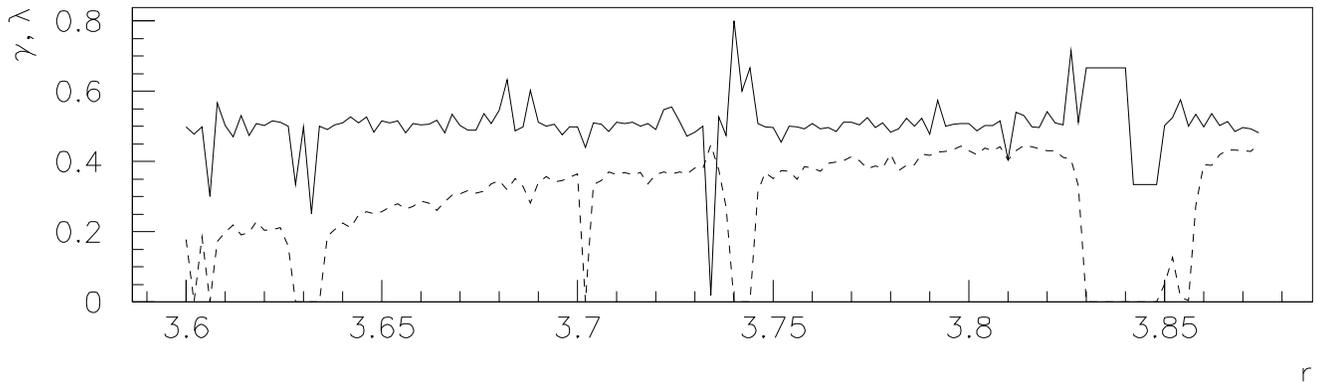
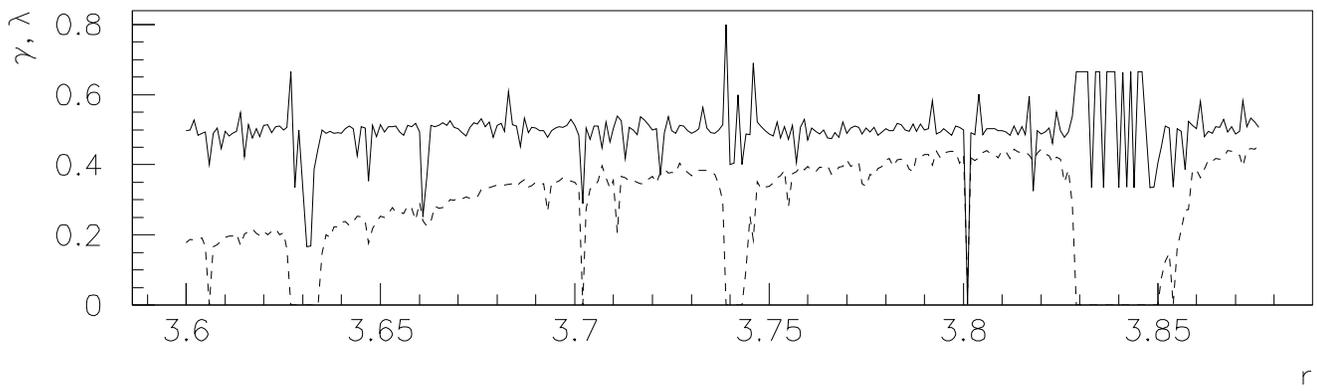
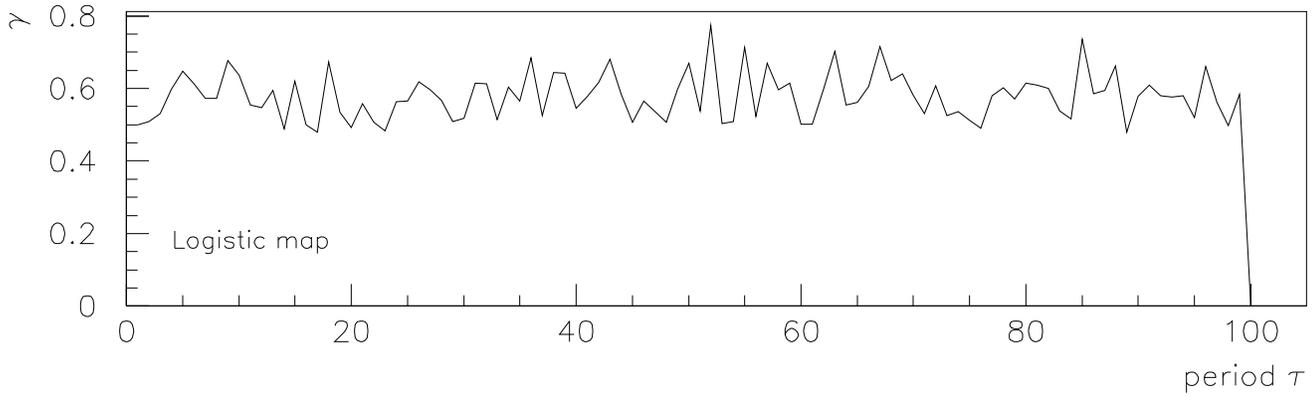
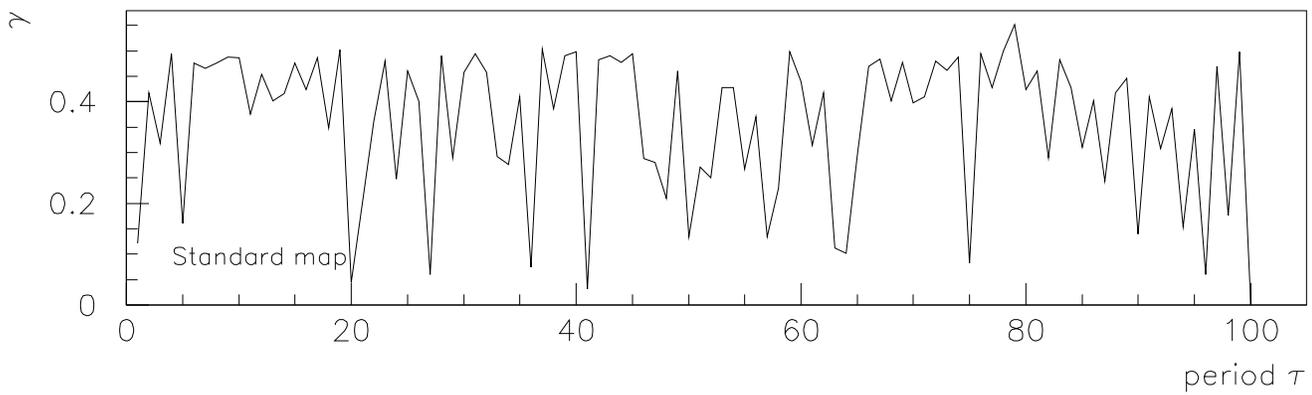
\end{document}